\documentclass[twocolumn,twoside]{IEEEtran}
\usepackage{amsmath,amssymb,amsthm,epsfig,color,subfigure,empheq,graphicx,pgfplots}
\usepackage{enumerate,url,algpseudocode,algorithm,wasysym,epstopdf,balance}

\newcommand \ba{\mathbf{a}}
\newcommand \bb{\mathbf{b}}

\newcommand \bn{\mathbf{n}}
\newcommand \br{\mathbf{r}}
\newcommand \bs{\mathbf{s}}

\newcommand \bv{\mathbf{v}}

\newcommand \bx{\mathbf{x}}

\newcommand \bz{\mathbf{z}}
\newcommand \bp{\mathbf{p}}
\newcommand \bq{\mathbf{q}}
\newcommand \bA{\mathbf{A}}
\newcommand \bB{\mathbf{B}}

\newcommand \bG{\mathbf{G}}
\newcommand \bSigma{\mathbf{\Sigma}}

\newcommand \hbSigma{\hat{\mathbf{\Sigma}}}
\newcommand \tbSigma{\tilde{\mathbf{\Sigma}}}

\newcommand \bI{\mathbf{I}}
\newcommand \bL{\mathbf{L}}

\newcommand \bF{\mathbf{F}}
\newcommand \bX{\mathbf{X}}
\newcommand \bY{\mathbf{Y}}
\newcommand \bZ{\mathbf{Z}}
\newcommand \bR{\mathbf{R}}
\newcommand \bW{\mathbf{W}}
\newcommand \bS{\mathbf{S}}
\newcommand \tbb{\tilde{\mathbf{b}}}
\newcommand \tbv{\tilde{\mathbf{v}}}

\newcommand \tbp{\tilde{\mathbf{p}}}
\newcommand \tbq{\tilde{\mathbf{q}}}
\newcommand \cbb{\check{\mathbf{b}}}

\newcommand \tbR{\tilde{\mathbf{R}}}
\newcommand \tbX{\tilde{\mathbf{X}}}

\newcommand \tbA{\tilde{\mathbf{A}}}

\newcommand \hbb{\hat{\mathbf{b}}}
\newcommand \mcL{\mathcal{L}}
\newcommand \mcB{\mathcal{B}}
\newcommand \bbeta{\boldsymbol{\beta}}

\DeclareMathOperator{\nullspace}{null}

\DeclareMathOperator{\diag}{dg}

\DeclareMathOperator{\trace}{Tr}

\DeclareMathOperator{\real}{Re}

\newtheorem{lemma}{Lemma}

\newtheorem{definition}{Definition}
\newtheorem{remark}{Remark}

\begin{document}
\title{Voltage Analytics for Power Distribution Network Topology Verification}
	
\author{
Guido Cavraro,
Vassilis Kekatos,~\IEEEmembership{Senior Member,~IEEE}, 
and 
Sriharsha Veeramachaneni


\thanks{G. Cavraro and V. Kekatos are with the Bradley Dept. of Electrical and Computer Engnr., Virginia Tech, Blacksburg, VA 24061, USA. S. Veeramachaneni is with Windlogics Inc., Saint Paul, MN 55108, USA. Emails: \{cavraro,kekatos\}@vt.edu, Sriharsha.Veeramachaneni@windlogics.com.}

}	


\maketitle

\begin{abstract}
Distribution grids constitute complex networks of lines oftentimes reconfigured to minimize losses, balance loads, alleviate faults, or for maintenance purposes. Topology monitoring becomes a critical task for optimal grid scheduling. While synchrophasor installations are limited in low-voltage grids, utilities have an abundance of smart meter data at their disposal. In this context, a statistical learning framework is put forth for verifying \textcolor{black}{single-phase} grid structures using non-synchronized voltage data. The related maximum likelihood task boils down to minimizing a non-convex function over a non-convex set. The function involves the sample voltage covariance matrix and the feasible set is relaxed to its convex hull. Asymptotically in the number of data, the actual topology yields the global minimizer of the original and the relaxed problems. Under a simplified data model, the function turns out to be convex, thus offering optimality guarantees. Prior information on line statuses is also incorporated via a maximum a-posteriori approach. The formulated tasks are tackled using solvers having complementary strengths. Numerical tests using real data on benchmark feeders demonstrate that reliable topology estimates can be acquired even with a few smart meter data, while the non-convex schemes exhibit superior line verification performance at the expense of additional computational time.
\end{abstract}


\section{Introduction}\label{sec:intro}
Low commercial interest together with the sheer coverage of residential low-voltage grids have resulted in their limited instrumentation. Traditionally, utility operators collect voltage, current, and power readings only from a few grid points. With the growing interest in integrating solar generation along with demand-response and electric vehicle charging programs, utilities need more refined models of their assets to accomplish grid scheduling tasks. To this end and given the currently prohibitive cost of installing synchrophasors on a wide scale, data from advanced metering infrastructure (AMI) and power inverters can provide useful grid information.

One such critical piece of information is the operational structure of a grid. Power networks are built with line redundancy for efficiency, reliability, and maintenance purposes. Although operators know the available lines along with their characteristics, the energized topology under which the grid operates at any given time may not be precisely known. At the transmission level, the grid topology is typically acquired using historical data, measurements, and the generalized state estimator~\cite{ExpConCanBook}. Detecting sudden single- and double-line outages via efficient enumerations has been suggested in~\cite{tate09pesgm}, while outages of multiple lines are unveiled via the sparse overcomplete representation of~\cite{ZhGi12}. Given power injections across the network over multiple times, grid connectivity is recovered via a blind matrix factorization in~\cite{LiPoSc13}. A Gaussian Markov random field has been postulated over bus voltage angles to localize transmission grid faults~\cite{HeZhang2011}. Instead of using electrical quantities, power system topologies are tracked using publicly available electricity prices in~\cite{KGB16}.

Focusing on distribution grids, reference \cite{ErTpVi13} exploits the time delays of power line communication signals to reveal the grid structure. The topology recovery scheme of \cite{BoSch13} relies on the properties of the inverse covariance matrix of bus voltage magnitudes. After developing a linearized distribution grid model, \cite{Deka1} generalize the previous schemes to voltage data from multiple feeders, correlated power injections, and grids with variable resistance-to-reactance ratios. The graph recovery algorithms of \cite{Deka1} have been extended to incorporate covariances of power injections from terminal nodes~\cite{Deka3}. Grid topology recovery has been tackled using graphical models by fitting a spanning tree based on the mutual information of voltage data~\cite{WengLiaoRajagopal17}. The aforementioned approaches rely on the \emph{ensemble} rather than \emph{sample} moments of meter data. In a grid of $N$ buses, the sample covariance of voltage data becomes invertible after collecting at least $N$ data. If meters report every few minutes, the sample covariance matrix can be inverted only after some hours. Even then, its inverse would deviate substantially from its ensemble counterpart.

Topology inference methods relying on synchrophasor data have also been suggested for distribution grids. The scheme of \cite{Arghandeh15} selects the topology attaining the best least-squares fit in an exhaustive fashion whose complexity grows exponentially in the number of configurations. A data-driven algorithm for detecting switching events based on topology signatures has been reported~\cite{DDA}. A sparse linear model capturing the voltage dependence between every node and all other nodes is sought via $\ell_1$-penalized regression in~\cite{LiaoWengRajagopal16Lasso}; yet the per-node models may not agree. In \cite{Ardakanian17}, the bus admittance matrix is found via linear regression and its observability is characterized presuming all non-metered buses do not inject power. Albeit synchrophasors for distribution grids are underway, their current cost inhibits wide adoption. 

The task of verifying grid topologies using non-synchronized voltage data is considered here. Our contribution is three-fold: First, topology verification in \textcolor{black}{single-phase} grids is posed as a maximum likelihood (ML) problem involving the \emph{sample} covariance matrix of voltage data. After reviewing the grid model in Section~\ref{sec:model}, the grid topology is captured by a binary vector (Section~\ref{subsec:detailed}). The associated non-convex set is relaxed to its convex hull and a stationary point of the non-convex likelihood function is found via gradient projection. Asymptotically in the number of data, the true topology constitutes the global minimizer for both problems. Second, the novel learning schemes simplify if lines are assumed to exhibit identical resistance-to-reactance ratios, see also \cite{BoSch13}. By further ignoring noise, the likelihood function becomes convex, and hence, numerical bounds on the suboptimality of the relaxation are obtained (Section~\ref{subsec:simplified}). Lastly, possible prior information on the status of individual lines is incorporated in terms of a maximum \emph{a posteriori} (MAP) estimator (Section~\ref{subsec:priors}). The numerical tests of Section~\ref{sec:tests} using actual data on benchmark feeders corroborate our findings, and conclusions are drawn in Section~\ref{sec:conclusions}.

Regarding \emph{notation}, lower- (upper-) case boldface letters denote column vectors (matrices), with the exception of the line complex power flow vector $\mathbf{S}$. Calligraphic symbols are reserved for sets. Symbol $^{\top}$ stands for transposition and \textcolor{black}{$\bar{z}$ is the complex conjugate of $z$}. Vectors $\mathbf{0}$ and $\mathbf{1}$ are the all-zero and all-one vectors. Symbol $\|\mathbf{x}\|_2$ denotes the $\ell_2$-norm of $\mathbf{x}$ and $\diag(\mathbf{x})$ defines a diagonal matrix having $\mathbf{x}$ on its diagonal. A symmetric positive (semi)definite matrix is denoted as $\mathbf{X}\succ \mathbf{0}$ $(\mathbf{X}\succeq \mathbf{0})$.

\section{Modeling Preliminaries}\label{sec:model}
The proposed topology verification schemes build on the approximate linearized model for distribution grids that are briefly reviewed next.
\allowdisplaybreaks
A radial single-phase grid having $N+1$ buses can be represented by a tree graph $\mathcal{T}=(\mathcal{N}_o,\mcL)$, whose nodes $\mathcal{N}_o:=\{0,\ldots,N\}$ correspond to buses and its edges $\mcL$ to distribution lines. The tree is rooted at the substation bus indexed by $n=0$, while every non-root bus $n\in\mathcal{N}:=\{1,\ldots,N\}$ is connected to its unique parent bus $\pi_n$ via line $n$. The grid is modeled by the branch flow equations~\cite{BW3}
\begin{subequations}\label{eq:msmv}
	\begin{align}
	s_n&=\sum_{k\in\mathcal{C}_n}S_k  - S_n + I_n z_n \label{eq:ms}\\
	v_n&=v_{\pi_n} - 2\real[\bar{z}_n S_n] +I_n |z_n|^2\label{eq:mv}\\
	|S_n|^2&=v_{\pi_n} I_n\label{eq:ml}
	\end{align}
\end{subequations}
where $z_n=r_n+jx_n$ is the impedance of line $n$; \textcolor{black}{$I_n$ is the squared current magnitude on line $n$}; $S_n$ is the complex power flow sent on line $n$ from the parent node $\pi_n$; and $\mathcal{C}_n$ is the set of children buses for $n$. Further, $v_n$ and $s_n=p_n+jq_n$ are the squared voltage magnitude and power injection at bus $n$. \textcolor{black}{Equation \eqref{eq:ms} follows from conservation of power; equation \eqref{eq:mv} describes the drop in squared voltages along line $n$~\cite{BW3}; and \eqref{eq:ml} stems from the definition of complex power flow.}

Since the nonlinearity in \eqref{eq:ml} complicates power flow calculations, distribution grids are sometimes studied using the \emph{linear distribution flow (LDF)} model~\cite{BW3}. \textcolor{black}{The latter follows from \eqref{eq:msmv} upon dropping the last summands in the right-hand sides (RHS) of \eqref{eq:ms}--\eqref{eq:mv}.} It can also be derived assuming voltage magnitudes to be close to unity and voltage angle differences between neighboring buses to be small~\cite{Deka1}; \textcolor{black}{or alternatively, via a first-order Taylor series approximation of power injections as functions of voltages~\cite{BoDo15}.}

\textcolor{black}{The topology of a grid with $N+1$ buses and $L$ distribution lines is captured by the branch-bus incidence matrix $\tbA \in\{0,\pm1\}^{L\times (N+1)}$ that can be partitioned into its first and the rest of its columns as $\tbA=[\mathbf{a}_0~\bA]$. For a radial grid $(L=N)$}, the \emph{reduced incidence matrix}~$\bA$ is square and invertible~\cite{GodsilRoyle}. \textcolor{black}{After ignoring losses, equations \eqref{eq:ms}--\eqref{eq:mv} can be compactly expressed in matrix-vector form as~\cite{VKZG16}}
\begin{subequations}\label{eq:mc}
	\begin{align}
	 \bs&=\bA^\top\bS\label{eq:mcs}\\
	 \bA\bv&=2\real[\diag(\bar{\bz}) \bS]  - \ba_0 v_0\label{eq:mcv}
	\end{align}
\end{subequations}
where $\bs=\bp+j\bq$ is the vector of complex nodal injections; $\bS$ is the vector of complex power flows; $\bz=\br+j\bx$ is the vector of line impedances; and $v_0$ is the squared voltage magnitude at the substation. \textcolor{black}{Substituting $\bS=\bA^{-\top}\bs$ from \eqref{eq:mcs} in \eqref{eq:mcv}, and exploiting the fact that $\mathbf{a}_0 +\bA\mathbf{1}=\mathbf{0}$ (due to $\tbA\mathbf{1}=\mathbf{0}$), \textcolor{black}{the vector $\bv$ collecting the squared voltage magnitudes at all buses in $\mathcal{N}$} can be approximated as~\cite{FCL13,Deka1,BoDo15}}
\begin{equation}\label{eq:model}
\bv\simeq \bR\bp + \bX\bq + v_0\mathbf{1}_N
\end{equation}
where the $N\times N$ matrices $\bR$ and $\bX$ are defined as
\begin{subequations}\label{eq:RX}
\begin{align}
\bR&:=2(\bA^\top \diag^{-1}(\br)\bA)^{-1}\label{eq:RX:R}\\
\bX&:=2(\bA^\top \diag^{-1}(\bx)\bA)^{-1}\label{eq:RX:X}.
\end{align}
\end{subequations}
The RHS of \eqref{eq:model} has been shown to underestimate the actual squared voltage magnitudes with the bias depending on $I_n$'s; see~\cite{GLTL12}. \textcolor{black}{Numerical tests though report approximation errors in voltage magnitudes less than 0.001 pu~\cite{GLTL12}. The approximate model of \eqref{eq:model} has been extended to multi-phase grids in~\cite{VKZG16}.}

Upon applying \eqref{eq:model} over two successive time instances, the fluctuations in squared voltage magnitudes $\tbv_t:=\bv_t-\bv_{t-1}$ caused by perturbations in power injection $\tbp_t:=\bp_t-\bp_{t-1}$ and $\tbq_t:=\bq_t-\bq_{t-1}$ follow the model
\begin{equation}\label{eq:modeld}
\tbv_t = \bR\tbp_t + \bX\tbq_t + \bn_t
\end{equation}
where the vector $\bn_t$ captures measurement noise, the approximation error introduced by LDF, and modeling inaccuracies.

\section{Topology Verification}\label{sec:verification}
\begin{figure}[t]
\centering
\includegraphics[width=0.40\textwidth]{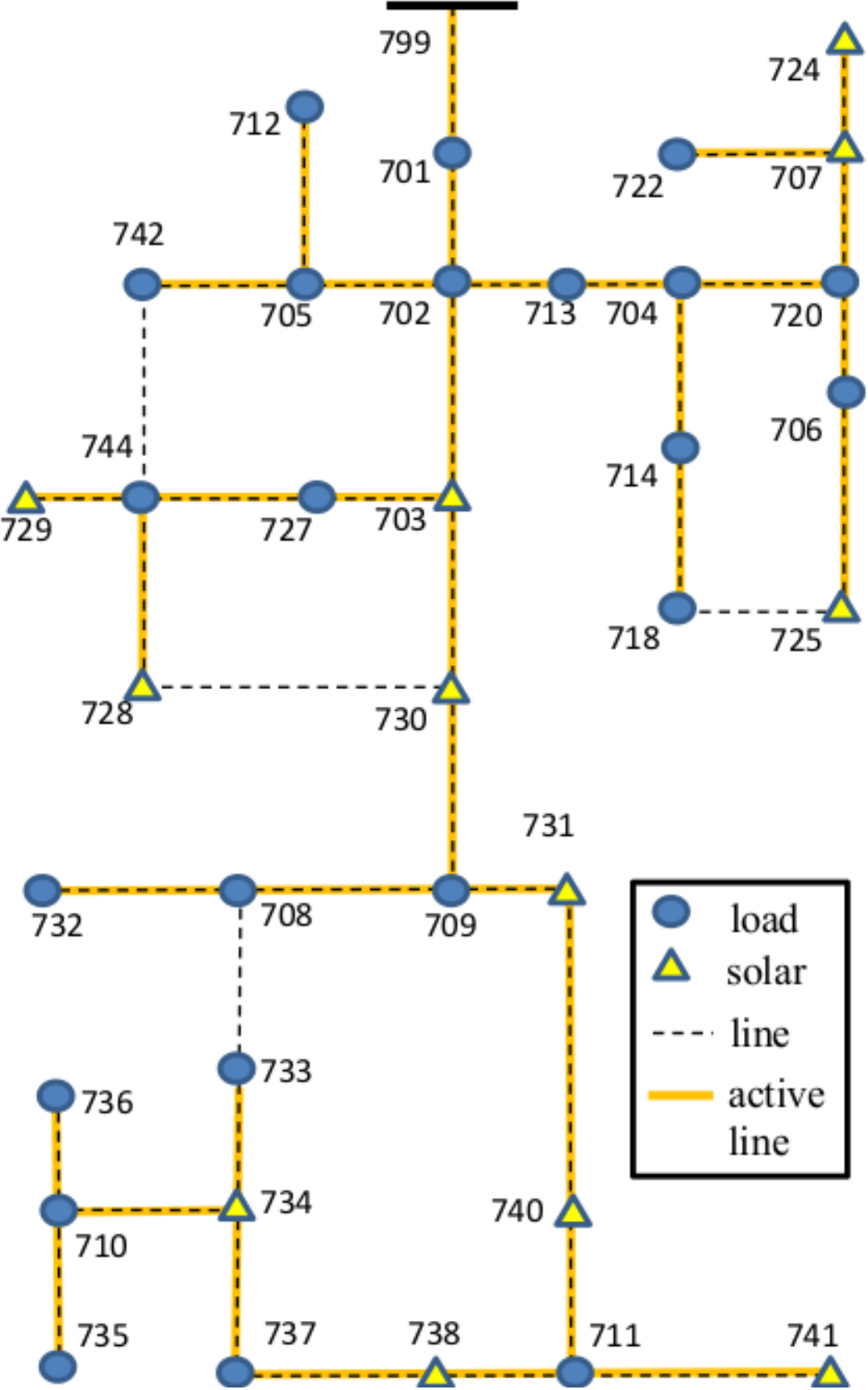}
\caption{Line infrastructure for the IEEE 37-bus feeder benchmark.}\label{fig:ieee37}
\end{figure}

Power grids are built with redundancy in line infrastructure. This redundancy improves system reliability against failures or during scheduled maintenance, while grids are regularly reconfigured for loss minimization~\cite{BW3}. Therefore, the set of energized lines $\mcL$ is a subset of the existing lines denoted by the set $\mcL_e$ with $|\mcL_e|=L_e$, and $\mcL \subseteq \mcL_e$. For example, line 742-744 in Fig.~\ref{fig:ieee37} belongs to $\mcL_e$ but not to $\mcL$.

To verify grid topologies using smart meter data, let us express $(\bR,\bX)$ defined in \eqref{eq:RX} in terms of the active lines as
\begin{equation}\label{eq:RX2}
\bR=2\left(\sum_{\ell \in\mcL} \frac{1}{r_\ell}\ba_\ell \ba_\ell^\top\right)^{-1}\text{and}~~ \bX=2\left(\sum_{\ell\in\mcL} \frac{1}{x_\ell}\ba_\ell \ba_\ell^\top\right)^{-1}
\end{equation}
where $\ba_\ell^\top$ is the $\ell$-th row of $\bA$. By slightly abusing notation, matrix $\bA$ here has been augmented to include both active and non-energized lines. \textcolor{black}{Let us introduce the variable $b_\ell$ for line $\ell$ taking the values $b_\ell=1$ if $\ell\in\mcL$; and $b_\ell=0$ otherwise, for $\ell\in\mcL_e$.} Collecting $b_\ell$'s in the binary $L_e$--length vector $\bb$, matrices $(\bR,\bX)$ can be written as
\begin{equation}\label{eq:RX3}
\bR^{-1}(\bb) =\frac{1}{2}\sum_{\ell \in\mcL_e} \frac{b_\ell}{r_\ell}\ba_\ell \ba_\ell^\top,~\bX^{-1}(\bb) =\frac{1}{2}\sum_{\ell\in\mcL_e} \frac{b_\ell}{x_\ell}\ba_\ell \ba_\ell^\top
\end{equation}
where now the summations extend over all lines in $\mcL_e$. Building on \eqref{eq:RX3}, verifying the grid topology entails finding $\bb$ from smart meter data. Meter readings are collected by the utility over $T$ time slots. Given the current and forthcoming communication capabilities of the smart metering infrastructure, samples are collected on a 5- to 15-min basis. Hence, within a period of one hour over which the grid topology is presumed unchanged, at most $T=4$ to $12$ voltage vectors can be recorded. 

Smart meter data comprise voltage magnitudes, active powers, and power factors among other quantities. Nonetheless, a smart meter monitors a single household, which may not necessarily correspond to a bus. The nodes of low-voltage grids are typically mapped to pole transformers, each one serving 5 to 10 residential costumers. Hence, power readings from meters may not be readily used. Instead, we statistically characterize power injections and rely on the second-order statistics of voltage readings to verify the grid topology.

\subsection{Verification under a detailed data model}\label{subsec:detailed}
\textcolor{black}{To arrive at a tractable statistical formulation, differential injection data $(\tbp_t,\tbq_t)$ are assumed to be zero-mean random vectors with covariance matrices $\bSigma_p:=\mathbb{E}[\tbp_t\tbp_t^\top]$, $\bSigma_q:=\mathbb{E}[\tbq_t\tbq_t^\top]$, and $\bSigma_{pq}:=\mathbb{E}[\tbp_t\tbq_t^\top]$. The noise vector $\bn_t$ is modeled as zero-mean independent Gaussian with variance $\sigma_n^2$, that is $\bn_t\sim \mathcal{N}(\mathbf{0},\sigma_n^2\bI)$. Under these assumptions and from \eqref{eq:modeld}, differential voltage data are also zero-mean with ensemble covariance matrix explicitly parameterized by the line indicator vector $\bb$ as}
\begin{align}
\bSigma(\bb) &=\bR(\bb)\bSigma_p \bR(\bb) + \bX(\bb)\bSigma_q \bX(\bb) \notag \\
&~~~+ \bR(\bb)\bSigma_{pq} \bX(\bb) + \bX(\bb)\bSigma_{pq}^\top \bR(\bb) + \sigma_n^2\bI_N \label{eq:C}.
\end{align}
Based on the preceding modeling, the topology verification problem is stated next.

\begin{definition}\label{def:tv}
The task of topology verification amounts to finding the subset $\mcL$ of active distribution lines given:\\
\hspace*{1em}(i) the distribution line infrastructure $\{\ba_\ell\}_{\ell\in\mcL_e}$;\\
\hspace*{1em}(ii) the line parameters $\{r_\ell,x_\ell\}_{\ell\in \mcL_e}$;\\
\hspace*{1em}(iii) the covariance matrices $\bSigma_p$, $\bSigma_q$, $\bSigma_{pq}$ ; and\\ 
\hspace*{1em}(iv) voltage readings $\{\tbv_t\}_{t=1}^T$ over $T$ time instances.
\end{definition}

\textcolor{black}{Different from the task of topology identification where both the existing lines and their impedances are unknown~\cite{BoDo15}, \cite{Deka1}, in topology verification the utility operator knows the existing line infrastructure. Here, the operator collects voltage readings, but power injections are described only through their first- and second-order moments. This is to cater the scenario where smart meter data from one house can be used to trace back the voltage magnitude at the related pole transformer. The net injection on the pole transformer serving 5-10 additional houses however cannot be inferred. The covariances can be estimated from historical data. Postulating structure on covariances (e.g., diagonal matrices) reduces the number of parameters.}

\textcolor{black}{To build a statistical model on differential voltage data, we invoke contemporary variants of the central limit theorem (CLT) and postulate that the probability density function (pdf) of each random vector $\tbv_t$ converges asymptotically in $N$ to a multivariate Gaussian pdf, even if $(\tbp_t,\tbq_t)$ are not Gaussian.}

\begin{remark}\label{re:CLT}
Note that $\tbv_t=\sum_{n=1}^N(\tilde{p}_{t,n} \br_n + \tilde{q}_{t,n}\bx_n)$, where $(\tilde{p}_{t,n},\tilde{q}_{t,n})$ are the $n$-th entries of $(\tbp_t,\tbq_t)$ and $(\br_n,\bx_n)$ are the $n$-th columns of $(\bR,\bX)$, accordingly. The $N$-dimensional vector $\tbv_t$ is therefore expressed as the sum of $N$ independent random vectors. Different from the classical multivariate CLT, the dimension of $\tbv_t$ grows with the number of its summands. Nevertheless, its pdf can still be approximated by a Gaussian pdf for bounded power injections~\cite[Sec.~2]{Chernozhukov}.
\end{remark}

\textcolor{black}{The statistical tests based on actual data in Sec.~\ref{sec:tests} validate this assumption. Thus, the pdf of $\tbv_t$ can be approximated as}
\begin{equation*}
p(\tbv_t;\bb)=\frac{|\bSigma(\bb)|^{-1/2}}{(2\pi)^{N/2}}\exp\left(-\frac{1}{2}\tbv_t^\top\bSigma^{-1}(\bb)\tbv_t \right).
\end{equation*}
\textcolor{black}{To fully characterize the collected data $\{\tbv_t\}_{t=1}^T$, their joint pdf should be provided. To that end, assuming independence across time is convenient and computationally tractable. In fact, the normalized autocorrelation functions depicted in Fig.~\ref{fig:xcorr} demonstrate that voltage data are relatively uncorrelated across time, especially for sampling periods larger than $T_s=5$~min. Due to Gaussianity, uncorrelatedness implies independence.} Therefore, the joint pdf for the entire voltage dataset becomes $p(\{\tbv_t\}_{t=1}^T;\bb)=\prod_{t=1}^T p(\tbv_t;\bb)$, or
\begin{equation*}
p(\{\tbv_t\}_{t=1}^T;\bb)=\frac{|\bSigma(\bb)|^{-T/2}}{(2\pi)^{NT/2}}\exp\left(-\frac{1}{2}\sum_{t=1}^T\tbv_t^\top\bSigma^{-1}(\bb)\tbv_t \right).
\end{equation*}
Upon observing the data $\{\tbv_t\}_{t=1}^T$, function $p(\{\tbv_t\}_{t=1}^T;\bb)$ becomes the likelihood function of the unknown line indicator vector $\bb$. Adopting a maximum likelihood (ML) approach, the sought $\bb$ can be found as the maximizer of $p(\{\tbv_t\}_{t=1}^T;\bb)$, or equivalently, as the minimizer of $-\log p(\{\tbv_t\}_{t=1}^T;\bb)$. After scaling by $2/T$ and ignoring constant terms, the negative log-likelihood function is
\begin{equation*}
f(\bb):=\log |\bSigma(\bb)| + \frac{1}{T}\sum_{t=1}^T\tbv_t^\top\bSigma^{-1}(\bb)\tbv_t.
\end{equation*}
Exploiting the linearity of the trace operator $\trace(\cdot)$, the function $f(\bb)$ can be equivalently expressed as
\begin{equation}\label{eq:f}
f(\bb)=\log |\bSigma(\bb)|+\trace(\bSigma^{-1}(\bb)\hbSigma)
\end{equation}
where $|\cdot|$ denotes the matrix determinant and $\hbSigma$ is the sample covariance matrix of the voltage data
\begin{equation}\label{eq:Sigmahat}
\hbSigma:=\frac{1}{T} \sum_{t=1}^T \tbv_t \tbv_t^\top.
\end{equation}
Function $f(\bb)$ consists of two terms. The term $\trace(\bSigma^{-1}(\bb)\hbSigma)$ aggregates all the information from voltage readings and can be understood as a data-fitting term. The term $\log|\bSigma(\bb)|$ on the other hand acts as a regularizer guarding $\bSigma(\bb)$ within the positive definite matrix cone and away from singularity~\cite{KGB16}.

Vector $\bb$ can now be found by minimizing
\begin{equation}\label{eq:MLE}
f^\star:=\min_{\bb\in \mcB, \mathbf{1}^\top \bb=L}f(\bb)
\end{equation}
where $\mcB:=\{0,1\}^{L_e}$ and the constraint $\mathbf{1}^\top \bb=L$ fixes the number of active lines to $L$. In a grid with $N+1$ nodes, requiring $L=N$ in \eqref{eq:MLE} enforces a radial structure. To see this, note that the graph induced by any $N$ edges not forming a tree is not a connected graph. The corresponding reduced incidence matrix $\bA(\bb)$ is not full column-rank. Consequently, matrices $(\bR(\bb),\bX(\bb))$ become singular and the cost function $f(\bb)$ goes to infinity.

Solving \eqref{eq:MLE} is non-trivial due to the binary constraints. The set $\mcB$ can be surrogated by its convex hull $\mcB_r:=[0,1]^{L_e}$ and a lower bound on $f^\star$ can be found by solving
\begin{equation}\label{eq:MLEr}
f^\star_r:=\min_{\bb\in\mcB_r,\mathbf{1}^\top \bb=L}~f(\bb).
\end{equation}

The function $f(\bb)$ can be shown to be non-convex in general. A projected gradient descent (PGD) scheme with a sufficiently small step size $\mu>0$ is guaranteed to converge to a point $\cbb$ satisfying the necessary optimality condition for \eqref{eq:MLEr}, that is~\cite[Prop.~6.1.3]{Be15}
\begin{equation}\label{eq:nc}
(\nabla f(\cbb))^\top (\bb - \cbb)\geq 0~\text{for all}~\bb\in \mcB_r~\text{and}~\mathbf{1}^\top \bb=L.
\end{equation}
After initializing $\bb$ at some $\bb^0$, the $k$-th PGD iteration reads
\begin{equation}\label{eq:PGD}
\bb^{k+1}:=\arg\min_{\bb\in \mcB_r, \mathbf{1}^\top \bb=L}\|\bb-\bb^k+\mu \nabla f(\bb^k)\|_2^2.
\end{equation}
The projection in \eqref{eq:PGD} can be handled either by a generic linearly-constrained convex quadratic program; by the lambda iteration method~\cite[Sec.~5.2.4]{ExpConCanBook}; or by dual ascent upon dualizing the constraint $\mathbf{1}^\top \bb=L$. The gradient $\nabla f(\bb)$ is calculated using standard matrix differential rules as shown in the Appendix.

\begin{lemma}\label{le:gradient}
The $\ell$-th entry of the gradient of $f(\bb)$ is
{\color{black}
\begin{align}\label{eq:gradient}
\frac{\partial f(\bb)}{\partial b_{\ell}}&=-\frac{1}{r_\ell}\ba_\ell^\top\bR(\bb)\bSigma_p\bR(\bb)\bF(\bb)\bR(\bb) \ba_\ell\nonumber\\
&\quad-\frac{1}{x_\ell}\ba_\ell^\top\bX(\bb)\bSigma_q\bX(\bb)\bF(\bb)\bX(\bb) \ba_\ell \nonumber\\
&\quad-\frac{1}{r_\ell}\ba_\ell^\top\bR(\bb)\bSigma_{pq}\bX(\bb)\bF(\bb)\bR(\bb) \ba_\ell\nonumber\\
&\quad-\frac{1}{x_\ell}\ba_\ell^\top\bX(\bb)\bF(\bb)\bR(\bb)\bSigma_{pq}\bX(\bb) \ba_\ell
\end{align}}
where $\bF(\bb):=\bSigma^{-1}(\bb)-\bSigma^{-1}(\bb) \hbSigma \bSigma^{-1}(\bb)$.
\end{lemma}

Consistent with the properties of the maximum likelihood estimator (see \cite{kay93book}), the true line indicator vector minimizes \eqref{eq:MLE} and \eqref{eq:MLEr}, when the number of data $T$ grows to infinity; see the Appendix for a proof.

\begin{lemma}\label{le:stationary}
Let $\bb_o$ be the true line indicator vector. If the sample covariance $\hbSigma$ has converged to the ensemble covariance $\bSigma(\bb_o)$, then $\bb_o$ is a stationary point of $f(\bb)$ and global minimizer for \eqref{eq:MLE} and \eqref{eq:MLEr}.
\end{lemma}

The PGD estimate $\cbb$ may not lie in the original non-convex set $\mcB$. Then, a feasible vector $\hbb$ can be heuristically obtained by selecting the lines corresponding to the $L$ largest entries of $\cbb$; by finding the minimum spanning tree on a graph having $\cbb$ as edge weights; or \textcolor{black}{by choosing $\cbb$ as the mean of independent Bernoulli distributions used to draw a number of line configurations from which the feasible one attaining the smallest $f(\bb)$ is returned as $\hbb$}. Obviously, it holds that $f^\star\leq f(\hbb)$ and $f(\cbb)\leq f(\hbb)$. Although $\hbb$ yields reasonable verification performance (see Sec.~\ref{sec:tests}), vector $\cbb$ may not be a global minimizer of \eqref{eq:MLEr} since $f(\bb)$ is non-convex.

\subsection{Verification under a simplified data model}\label{subsec:simplified}
To alleviate the issues with the non-convexity of $f(\bb)$, we alternatively pursue a simplified model for voltage data by resorting to two assumptions~\cite{BoDo15}:\\
\hspace*{1em}(a1) The resistance-to-reactance ratios $\alpha_\ell:=r_\ell / x_\ell$ are identical for all lines, i.e., $\alpha_\ell=\alpha$ for all $\ell\in\mcL_e$.\\
\hspace*{1em}(a2) The noise term $\bn_t$ is negligible, or $\sigma_n^2=0$.\\
Regarding (a1), the ratios $\alpha_\ell$'s practically tend to lie within a limited range $[\alpha_{\min},\alpha_{\max}]$ as listed on Table~\ref{tbl:linedata}. The parameter $\alpha$ can be picked as the average or the median of the known actual ratios $\alpha_\ell$'s. Ignoring the noise term simplifies considerably the approach.

\begin{table}
\renewcommand{\arraystretch}{1}
\caption{Distribution line resistance-to-reactance ratios}\vspace*{-0.5em}
\label{tbl:linedata} \centering
\begin{tabular}{|lccccc|}
\hline
Feeder & $\alpha_{\min}$ & $\alpha_{\max}$ & mean & std & median \\
\hline
IEEE 34-bus & 1.00 & 1.88 & 1.41 & 0.29 & 1.37\\
IEEE 37-bus & 1.48 & 2.70 & 2.72 & 0.45 & 1.93\\
IEEE 123-bus & 0.42 & 2.02 & 0.74 & 0.38 & 0.97\\
\hline
\end{tabular}
\vspace*{-0.5em}
\end{table}

Under (a1)--(a2), voltage data are modeled as
\begin{equation}\label{eq:modela}
\tbv_t = \alpha\bX\tbp_t + \bX\tbq_t.
\end{equation}
Defining \textcolor{black}{$\bSigma_\alpha:=\alpha^2 \bSigma_p+\bSigma_q+\alpha(\bSigma_{pq}+\bSigma_{pq}^\top)$}, their ensemble covariance becomes
\begin{equation}\label{eq:Ca}
\tbSigma(\bb) = \bX(\bb)\bSigma_\alpha \bX(\bb).
\end{equation}
Using properties of the determinant and trace, the original negative log-likelihood in \eqref{eq:f} is surrogated by the function
\begin{equation}\label{eq:tf}
\tilde{f}(\bb):=-2\log |\bX^{-1}(\bb)|+\trace[\bX^{-1}(\bb)\bSigma_\alpha^{-1}\bX^{-1}(\bb)\hbSigma]
\end{equation} 
which enjoys convexity as shown in the Appendix.

\begin{lemma}\label{le:tfconvex}
The function $\tilde{f}(\bb)$ is strictly convex over $\bb\in \mathbb{R}^L_+$.
\end{lemma}

Based on this simplified voltage data model, the line indicator vector can be recovered as the minimizer of 
\begin{equation}\label{eq:MLEtilde}
\tilde{f}^\star:=\min_{\bb\in \mcB, \mathbf{1}^\top \bb=L}~\tilde{f}(\bb).
\end{equation}
Again, to deal with the non-convexity of the feasible set, function $\tilde{f}(\bb)$ is now minimized over its convex hull instead
\begin{equation}\label{eq:MLErc}
\tbb:=\arg\min_{\bb\in \mcB_r, \mathbf{1}^\top \bb=L}~\tilde{f}(\bb).
\end{equation}
Different from \eqref{eq:MLEr}, the optimization in \eqref{eq:MLErc} is convex. \textcolor{black}{Upon replacing $\bSigma(\bb)$, $\bR(\bb)$, and $r_\ell$, respectively by $\tbSigma(\bb)$, $\alpha \bX(\bb)$, and $\alpha x_\ell$ onto \eqref{eq:gradient}, the $\ell$-th entry of $\nabla \tilde{f}(\bb)$ becomes}
\begin{equation*}
\frac{\partial \tilde{f}(\bb)}{\partial b_{\ell}}=\frac{1}{x_\ell}\ba_\ell^\top[\bSigma_\alpha^{-1}\bX^{-1}(\bb)\hbSigma-\bX(\bb) ]\ba_\ell,\quad \forall \ell\in\mcL_e.
\end{equation*} 
It is easy to verify that under the model of \eqref{eq:modela}, Lemma~\ref{le:stationary} carries over to function $\tilde{f}$. Precisely, \textcolor{black}{if $\bb_o$ is the true line indicator vector} and $\hbSigma$ has converged to the ensemble covariance $\tbSigma(\bb_o)=\bX(\bb_o)\bSigma_\alpha \bX(\bb_o)$ under the approximate model, then $\nabla \tilde{f}(\bb_o)=\mathbf{0}$. The latter together with the strong convexity of $\tilde{f}$ imply that $\bb_o$ is the unique minimizer of \eqref{eq:MLEtilde} and \eqref{eq:MLErc}.

The minimizer of \eqref{eq:MLErc} can be found by the PGD iteration
\begin{equation}\label{eq:PGD2}
\bb_r^{k+1}=\arg\min_{\bb\in \mcB_r, \mathbf{1}^\top \bb=L}\|\bb-\bb^k+\mu \nabla \tilde{f}(\bb^k)\|_2^2.
\end{equation}
The Frank-Wolfe (FW) or conditional gradient scheme is a competing alternative for tackling \eqref{eq:MLErc}; see also \cite{ZKG16,jaggi2013revisiting}. The scheme is suitable for minimizing a differentiable convex function over a compact set. Tailoring the FW iterates to the problem at hand provides
\begin{subequations}\label{eq:FW}
\begin{align}
\breve{\bb}_r^{k+1}&\in\arg\min_{\bb\in \mathcal{B}_r,\mathbf{1}^\top \bb=L} \mathbf{b}^\top \nabla \tilde{f}(\bb_r^k)\label{eq:FWa}\\
\bb_r^{k+1}&=\bb_r^k +\tfrac{2}{k+2}(\breve{\bb}_r^{k+1} - \bb_r^k).\label{eq:FWb}
\end{align}
\end{subequations}
\textcolor{black}{The step in \eqref{eq:FWa} involves a linear program (LP) whose minimizer can be found by simply setting the entries of $\breve{\bb}_r^{k+1}$ corresponding to the $L$ smallest entries of $\nabla \tilde{f}(\bb_r^k)$ to one, and setting the remaining entries to zero. The FW scheme converges sublinearly in $\mathcal{O}(1/k)$ iterations, while the PGD scheme exhibits linear convergence~\cite{Be15}. Nonetheless, the LP step of FW in \eqref{eq:FWa} is computationally cheaper than the costly projection in \eqref{eq:PGD2}. Practically, the FW scheme is preferred over PGD when a low-accuracy solution suffices~\cite{jaggi2013revisiting}.}

If the minimizer $\tbb$ of \eqref{eq:MLErc} belongs to $\mcB$, it is a minimizer of \eqref{eq:MLEtilde} too. Otherwise, a suboptimal $\bb'$ can be obtained from $\tbb\in\mcB_r$ using the heuristics discussed earlier. In this case however, thanks to the convexity in \eqref{eq:MLErc}, the suboptimality gap can be bounded as $f(\tbb)\leq \tilde{f}^\star\leq f(\bb')$; see also~\cite[Sec.~7.5]{BoVa04}.

Comparing the two approaches, the one in Sec.~\ref{subsec:detailed} relies on the detailed model of \eqref{eq:modeld} and fully exploits the line data $\{r_\ell,x_\ell\}$. However, it entails a non-convex negative log-likelihood $f(\bb)$ that can be minimized only up to a point satisfying \eqref{eq:nc}. The approach of this section builds on the simplified model of \eqref{eq:modela} requiring only $\{x_\ell\}$ and the approximate common ratio $\alpha$. Because the associated function $\tilde{f}(\bb)$ is convex, the relaxation in \eqref{eq:MLErc} can be solved to optimality and with quantifiable gap from the non-convex task of \eqref{eq:MLEtilde}. 

\section{Extensions}\label{sec:extensions}
Our learning toolbox is extended next to cope with meshed topologies, disjoint feeders, and prior information on lines.

\subsection{Meshed topologies}\label{subsec:meshed}
Up to this point, the underlying grid topology was presumed to be radial. However, urban power distribution grids may exhibit operationally meshed structures. For distribution grids with $L\geq N$, the model in \eqref{eq:modeld} generalizes to~\cite{Deka1}, \cite{BoDo15}
\begin{equation}\label{eq:modelm}
\tbv_t = 2 (\bG + \bB\bG^{-1}\bB)^{-1}\tbp_t + 2 (\bB + \bG\bB^{-1}\bG)^{-1} \tbq_t + \textcolor{black}{\bn_t}
\end{equation}
for the bus conductance and susceptance matrices
\begin{subequations}\label{eq:BG}
\begin{align}
\bG&:=\bA^\top \diag\left(\left\{\frac{r_\ell}{r_\ell^2+x_\ell^2}\right\}_{\ell}\right)\bA\label{eq:BG:B}\\
\bB&:=\bA^\top \diag\left(\left\{\frac{x_\ell}{r_\ell^2+x_\ell^2}\right\}_{\ell}\right)\bA.\label{eq:BG:G}
\end{align}
\end{subequations}
It is not hard to show that for radial networks where $L=N$, matrix $\bA$ is square and the model in \eqref{eq:modelm} simplifies to \eqref{eq:modeld}.

\textcolor{black}{The scheme of Sec.~\ref{subsec:detailed} can be extended to meshed networks based on \eqref{eq:modelm}. The gradient of $f(\bb)$ can be derived by adopting the proof of Lemma~\ref{le:gradient} as outlined next. The formulas of \eqref{eq:diff1}--\eqref{eq:diff2} carry over as long as $\bR$ is substituted by $\tbR:= 2(\bG + \bB\bG^{-1}\bB)^{-1}$, and $\bX$ is substituted by $\tilde{\bX}:= 2(\bB + \bG\bB^{-1}\bG)^{-1}$. In this case however, the partial derivatives in \eqref{eq:diff3} are replaced by
\begin{align*}
\frac{\partial \tilde\bR(\bb)}{\partial b_{\ell}}& = -\frac{2r_\ell}{r_\ell^2 + x_\ell^2} \tbR(\bb)\ba_\ell \ba_\ell^\top\tbR(\bb)\\
&~~~+ \frac{2r_\ell}{r_\ell^2 + x_\ell^2}\bB(\bb)\bG^{-1}(\bb) \ba_\ell \ba_\ell^\top \bG^{-1}(\bb) \bB(\bb)\\
&~~~+ \frac{2x_\ell}{r_\ell^2 + x_\ell^2}\tbR(\bb)\bB(\bb)\bG^{-1}(\bb) \ba_\ell \ba_\ell^\top \tbR(\bb)\\
&~~~+\frac{2x_\ell}{r_\ell^2 + x_\ell^2}\tbR(\bb)\ba_\ell \ba_\ell^\top \bG^{-1}(\bb) \bB(\bb)\tbR(\bb)
\end{align*}
for all $\ell$. An analogous expression holds for $\frac{\partial \tbX(\bb)}{\partial b_{\ell}}$.}
 
On the other hand, under assumptions (a1)--(a2), the matrices in \eqref{eq:BG} simplify as $\bG=\frac{2\alpha}{\alpha^2+1}\bX^{-1}$ and $\bB=\frac{2}{\alpha^2+1}\bX^{-1}$. Simple algebraic manipulations show that \eqref{eq:modela}--\eqref{eq:tf} and the schemes of Sec.~\ref{subsec:simplified} apply \emph{unaltered} to meshed grids. This point further justifies the benefit of adopting the approximate model under (a1)--(a2).

\subsection{Disjoint feeders and multiple substation buses}\label{subsec:forest}
Voltage readings may be collected from multiple radial feeders or a meshed grid rooted at multiple substation buses~\cite{Deka1}. These scenarios can be handled in a uniform fashion under a simple modification. It is known that the nullity of the Laplacian matrix $\bL:=\tbA^\top \diag^{-1}(\br)\tbA$, i.e., the number of its zero eigenvalues, is equal to the number of the connected components in the underlying graph~\cite{GodsilRoyle}. Moreover, every principal minor of $\bL$ is strictly positive definite if the graph is connected. By definition, the nullity of $\bL$ coincides with the nullspace dimension of $\tbA$. 

When the grid is a single connected feeder, then $\dim(\nullspace(\tbA))=1$. Dropping any column of $\tbA$ provides a full column-rank matrix. Recall that $\bA$ is derived from $\tbA$ by removing the column corresponding to the substation bus.

The idea extends naturally. If data are collected from different feeders, the corresponding Laplacian can be transformed under column- and row-permutations to a block diagonal matrix whose diagonal blocks constitute the Laplacian matrices of the connected grid. To guarantee the invertibility of the individual Laplacians, the matrix $\bA$ used in our schemes is obtained from $\tbA$ by dropping the columns related to all substation buses. 

\subsection{Incorporating prior information on lines}\label{subsec:priors}
In meshed grids, the utility may not know the exact number of active lines $L$. \textcolor{black}{On the other hand, prior information on individual lines being active could be known through a generalized state estimator; current magnitude readings on transformers; or commonly used grid configurations (e.g., seasonally).} To cope with these setups, a maximum \emph{a-posteriori} (MAP) approach is adopted. The indicator $b_\ell$ for line $\ell$ is modeled as a Bernoulli random variable with given mean $\mathbb{E}[b_\ell]=\pi_\ell$. The prior probability mass function (pmf) for $b_\ell$ can be conveniently expressed as $\Pr(b_\ell)=\pi_\ell^{b_\ell}(1-\pi_\ell)^{1-b_\ell}$. To ease modeling, it is postulated that lines are active independently. Then, the joint pmf for $\bb$ is $\Pr(\bb)=\prod_{\ell\in\mcL_e}\Pr(b_\ell)$ up to normalization and its negative logarithm is
\begin{align}\label{eq:pmf}
-\log \Pr(\bb)&=-\sum_{\ell\in\mcL_e} b_\ell\log \pi_\ell +(1-b_\ell)\log(1-\pi_\ell)\nonumber\\
&=\sum_{\ell\in\mcL_e} b_\ell \beta_\ell -\log(1-\pi_\ell) 
\end{align}
where $\beta_\ell:=\log\left(\frac{1-\pi_\ell}{\pi_\ell}\right)$ for all $\ell\in\mcL_e$.

The MAP estimate for $\bb$ is defined as the vector maximizing the posterior pdf $p(\bb|\{\tbv_t\}_{t=1}^T)$. By Bayes' rule, the latter is proportional to the product $p(\{\tbv_t\}_{t=1}^T;\bb)\Pr(\bb)$ scaled by the inconsequential term $p(\{\tbv_t\}_{t=1}^T)$. The MAP estimate can be equivalently found by minimizing
\[-\log p(\bb|\{\tbv_t\}_{t=1}^T)=-\log p(\{\tbv_t\}_{t=1}^T;\bb) -\log \Pr(\bb).\]
By collecting the given parameters $\beta_\ell$'s in vector $\bbeta$ and ignoring constants, the latter leads to the problem
\begin{equation}\label{eq:MAP}
\bb_\text{MAP}:=\arg\min_{\bb\in \mcB} ~ \tfrac{T}{2}f(\bb) + \bbeta^\top \bb.
\end{equation}
If the simplified model of \eqref{eq:modela} is used in lieu of \eqref{eq:modeld}, the function $f(\bb)$ in \eqref{eq:MAP} is replaced by the convex $\tilde{f}(\bb)$. Contrasted to \eqref{eq:MLEtilde}, the task in \eqref{eq:MAP} leverages prior information on lines
\begin{itemize}
\item If line $\ell$ is likely to be active, then $\pi_\ell>1/2$ and $\beta_\ell<0$.
\item If line $\ell$ is known to be active, then $\pi_\ell=1$ and $\beta_\ell=-\infty$, thus forcing the $\ell$-th entry of $\bb_\text{MAP}$ to one.
\item No prior on line $\ell$ entails $\pi_\ell=1/2$ and $\beta_\ell=0$.
\item \textcolor{black}{If all lines statuses exhibit the same mean $\pi_c$, that is $\mathbb{E}[b_\ell]=\pi_\ell=\pi_c$ for $\ell\in\mcL_e$}, the linear term in the cost of \eqref{eq:MAP} becomes $\pi_c\mathbf{1}^\top \bb$. The latter can be interpreted as the Lagrangian counterpart of \eqref{eq:MLEtilde}.
\end{itemize}


The feasible set $\mcB$ can be again relaxed to $\mcB_r$. Lacking the coupling constraint $\mathbf{1}^\top\bb=L$, the PGD update
\begin{equation}\label{eq:PGD3}
\bb^{k+1}=\arg\min_{\bb\in \mcB_r}\|\bb-\bb^k+\mu \nabla f(\bb^k)+\mu\bbeta\|_2^2
\end{equation}
decouples across lines and the $\ell$-th entry of $\bb^{k+1}$ is found as
\begin{equation}\label{eq:PGD4}
b_\ell^{k+1}:=\left[ b_\ell^k -\mu\frac{\partial f(\bb^k)}{\partial b_\ell} -\mu \beta_\ell\right]_{0}^1
\end{equation}
where $[x]_0^1$ projects $x$ onto the interval $[0,1]$.

\section{Numerical Tests}\label{sec:tests}
\begin{figure}[t]
\centering	
\includegraphics[width=0.48\textwidth]{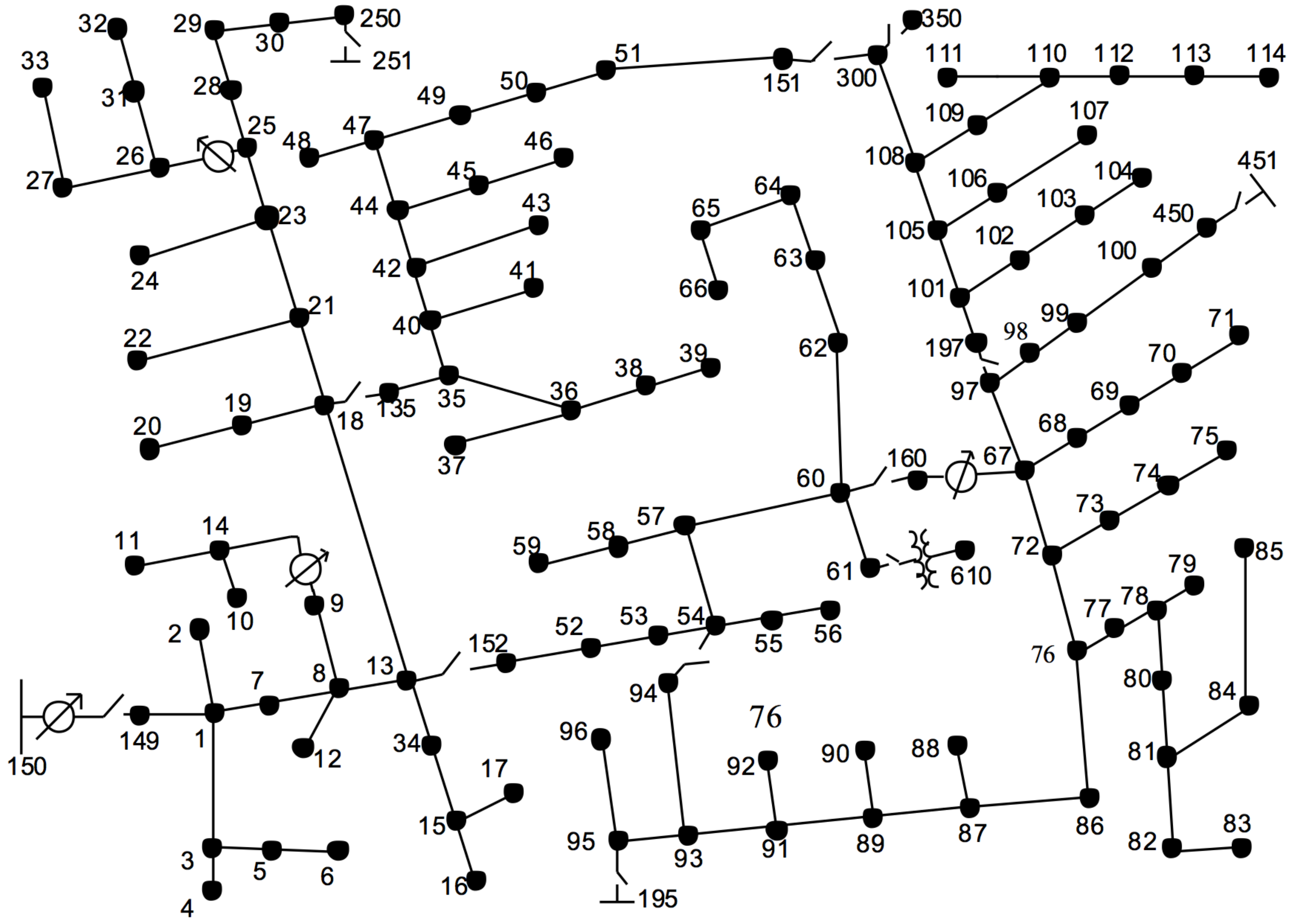}
\caption{The IEEE 123-bus distribution feeder with nine switches.}
\vspace*{-.5em}
\label{fig:ieee123}
\end{figure}

The developed schemes were validated using the IEEE 37-bus and 123-bus feeders depicted in Figures~\ref{fig:ieee37} and \ref{fig:ieee123}~\cite{Kersting}. \textcolor{black}{The feeders were modified to single-phase using the process described in~\cite{GLTL12}. Buses $\{150, 195,251,451\}$ served as possible substations for the IEEE 123-bus feeder. The actual topologies for both grids were randomly chosen allowing for both radial and meshed configurations. Each bus hosted a random number of four to ten houses, some of them hosting photovoltaic panels. The PV inverters implemented the following reactive control rule: if the generated active power exceeded half the PV capacity, the inverter would absorb reactive power according to an increasing function of its active power~\cite{Samadi2014}.}

Voltages were obtained by running a power flow solver using actual load and solar data from the Pecan Street dataset collected during January 1, 2013~\cite{pecandata}. Since this dataset does not provide reactive loads, power factors were drawn uniformly within 0.95 leading to 0.95 lagging. The measurement noise was modeled as zero-mean Gaussian with a 3-sigma deviation matching 0.5\% of the actual value. This is consistent with a study where 99\% of the tested meters attained relative accuracy of less than 0.5\%~\cite{an2011smart}. The schemes were run using MATLAB on a 2.7 GHz Intel Core i5 with 12GB RAM. 

\begin{figure}[t]
\centering	
\includegraphics[width=0.48\textwidth]{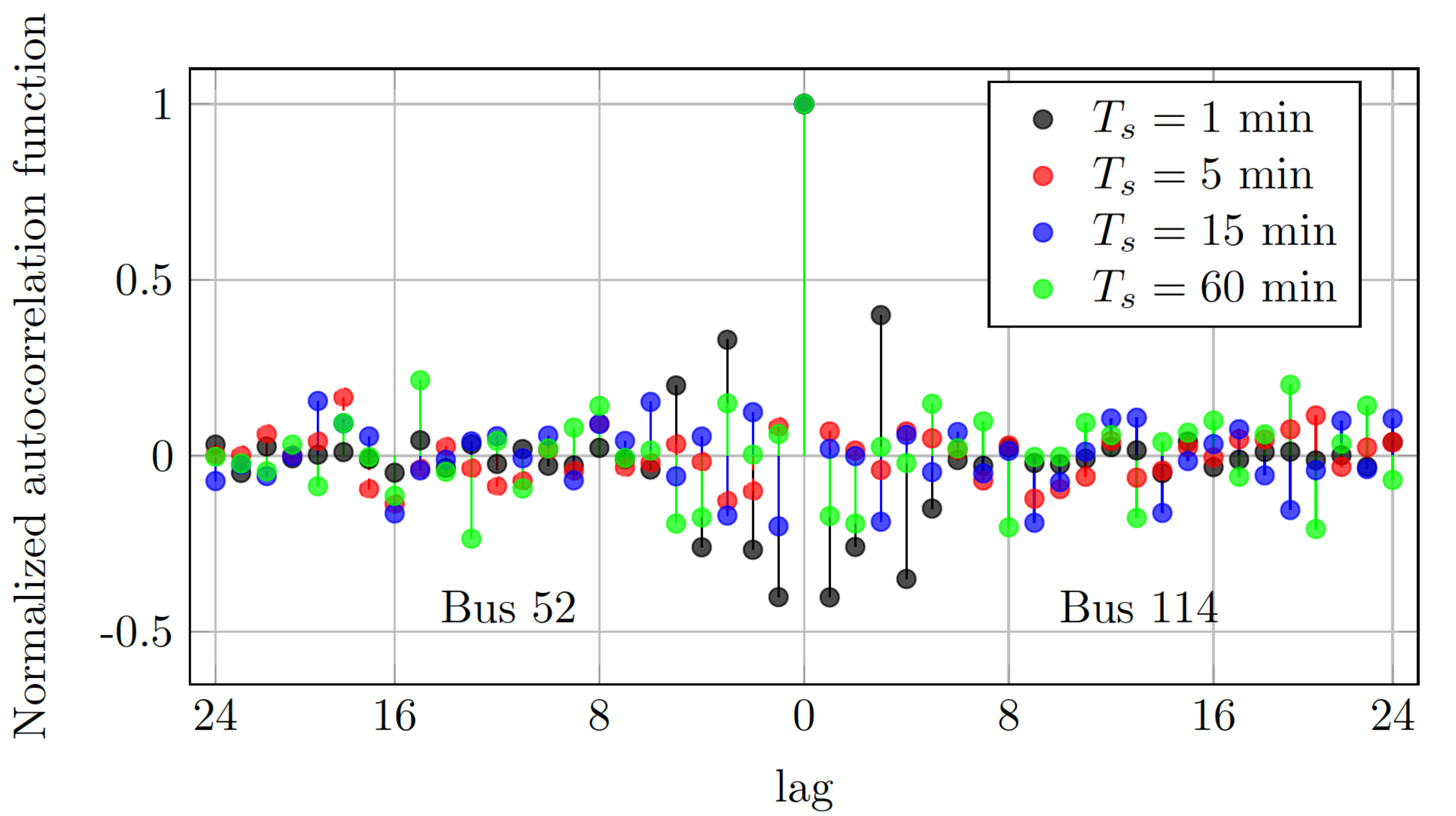}
\caption{Normalized autocorrelation function for voltage data on buses 52 (left) and 114 (right) computed using one-day data on the IEEE 123-bus feeder.}
\vspace*{-.5em}
\label{fig:xcorr}
\end{figure}

\textcolor{black}{First, the Gaussianity of voltage data was tested through the Kolmogorov-Smirnov test for various sampling rates using one-day data from the IEEE 123-bus feeder. Table~\ref{tbl:KStest} shows that the Gaussianity hypothesis is well suited when $T_s$ is 15, 30, and 60~min, but it is less accurate for $T_s=5$~min. Nevertheless, the numerical tests presented later demonstrate that the developed schemes exhibit reasonable performance even for $T_s=5$~min. Having Gaussianity in place, independence is ensured if data are uncorrelated over time. Figure~\ref{fig:xcorr} plots the normalized autocorrelation function for voltage data at buses 52 and 114 over $T_s$, which decays fast.}

\begin{table}
\renewcommand{\arraystretch}{1}
\caption{Number of buses on the IEEE 123-bus feeder that failed to pass the Kolmogorov-Smirnov test.}\vspace*{-0.5em}
\label{tbl:KStest} \centering
\begin{tabular}{|c|c|c|c|}
\hline
$T_s = 5$~min & $T_s = 15$~min & $T_s = 30$~min & $T_s = 60$~min \\
\hline
47 & 3 & 1 & 0 \\
\hline
\end{tabular}
\vspace*{-1em}
\end{table}

The ML verification schemes of Sec.~\ref{sec:verification} were then evaluated using 50 Monte Carlo runs over random topologies. The sampling rate was $T_s=5$~min. The covariances $\bSigma_p$, $\bSigma_q$, and $\bSigma_{pq}$ were approximated as diagonal matrices whose diagonal entries were estimated from historical data. In the first setup, line impedances were known and statuses were verified under the detailed model by solving \eqref{eq:MLEr}. In the second setup, line resistance-to-reactance ratios $\{\alpha_\ell\}_{\ell\in\mcL_e}$ were approximated by their average, and line verification was tackled via \eqref{eq:MLErc}. Because \eqref{eq:MLEr} is non-convex, the solution found by the PGD iterates of \eqref{eq:PGD} depends on the initialization $\bb^0$. The latter was selected as the minimizer of the convex problem in \eqref{eq:MLErc}. The solutions to \eqref{eq:MLEr} and \eqref{eq:MLErc} were projected onto $\mcB$ by setting their $L$ largest entries to one and nulling the rest.

\begin{figure}[t]
\centering	
\includegraphics[width=0.47\textwidth]{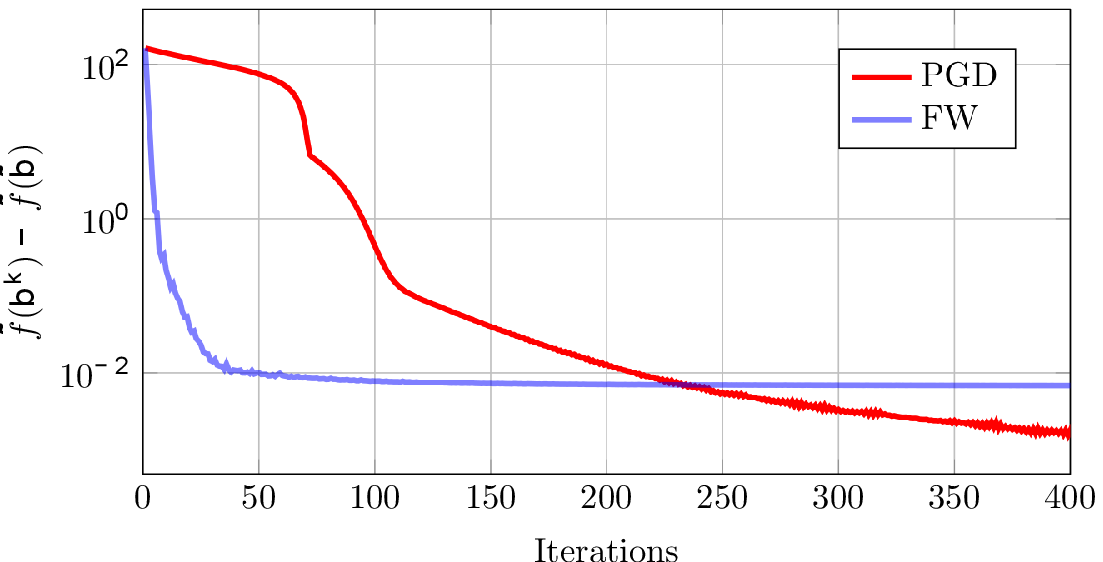}
\vspace*{1em}
\includegraphics[width=0.45\textwidth]{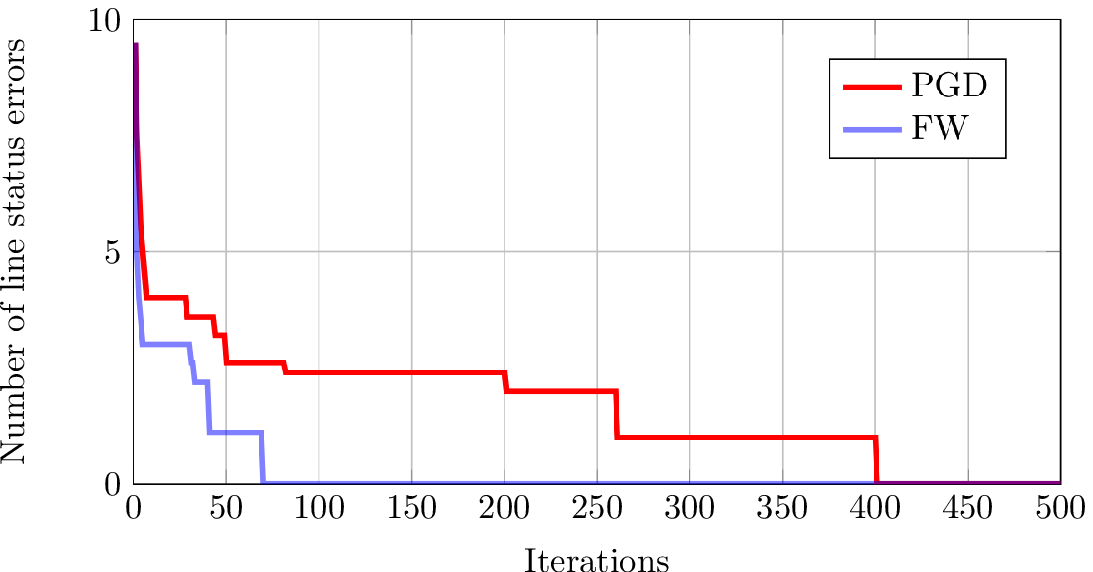}
\vspace*{-1em}
\caption{(top) Convergence of the solvers of \eqref{eq:MLErc} for the IEEE 37-bus grid; (bottom) Detection performance after projecting the intermediate iterates.}
\label{fig:thrs}
\end{figure}

\begin{figure}[t]
\centering	
\includegraphics[width=0.46\textwidth]{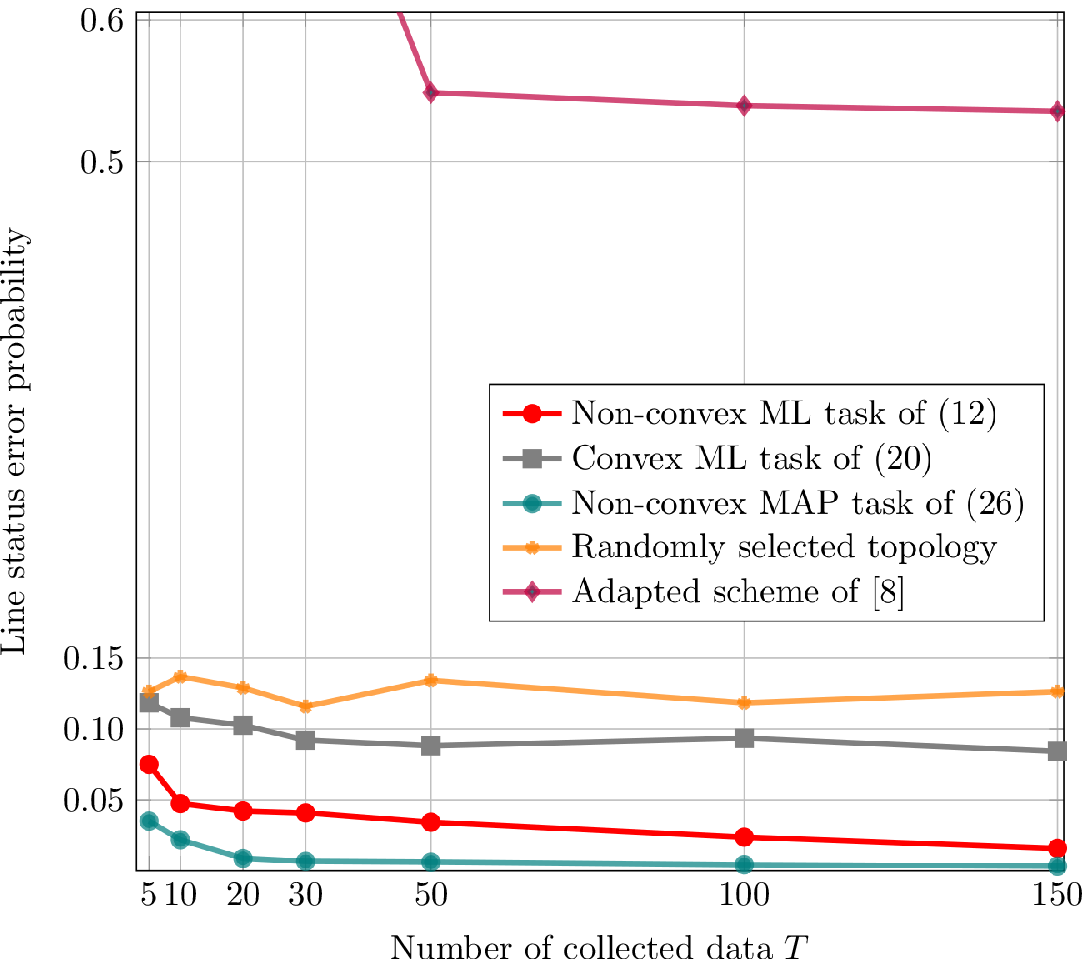}\\
\vspace*{.5em}
\includegraphics[width=0.46\textwidth]{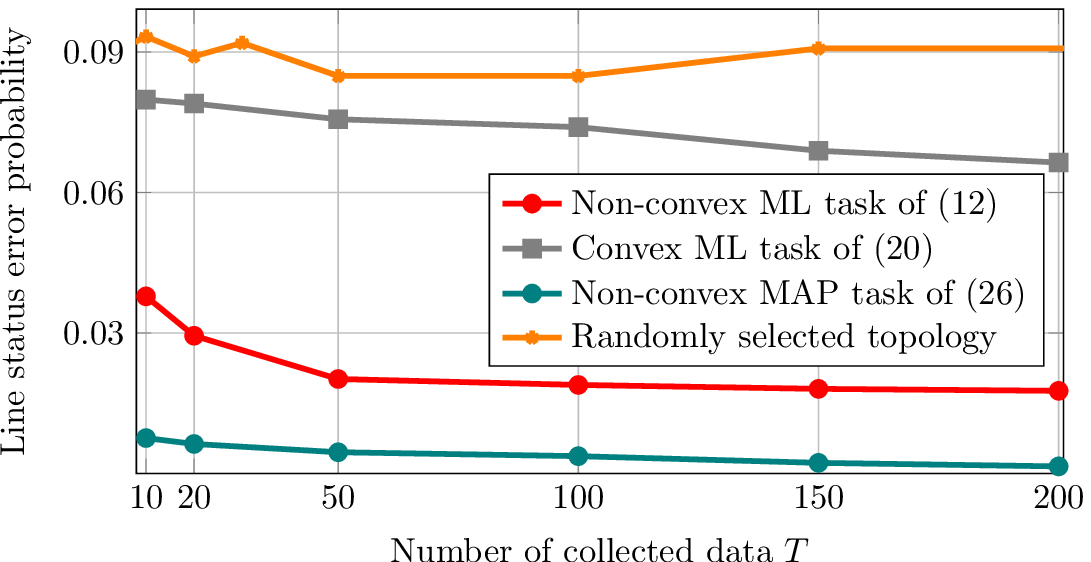}
\caption{Probability of line status errors for the IEEE 37-bus (top) and the IEEE 123-bus (bottom) feeders with $T_s = 5$~min.}
\label{fig:res}
\end{figure}

The convex problem in \eqref{eq:MLErc} was solved using both the PGD iterates of \eqref{eq:PGD}, and the FW iterates of \eqref{eq:FW}. The non-convex problem of \eqref{eq:MLEr} was solved solely by the PGD scheme, since the FW iterates are guaranteed to converge only under a non-trivial step size selection~\cite{Lacoste17}. The top panel of Fig.~\ref{fig:res} shows the convergence of the PGD and FW iterates for a particular instance of \eqref{eq:MLErc} with $T=100$. \textcolor{black}{Although solving the LP is computationally cheaper than the PGD projection, the execution time in both schemes is dominated by calculating $\nabla \tilde{f}(\bb)$. Numerical tests verify that the per-iteration time is similar for PGD and FW. Although PGD converges faster overall, the FW iterates exhibit faster initial convergence. The bottom panel of Fig.~\ref{fig:thrs} shows the number of line status errors, if the intermediate estimates $\bb_r^{k}$ are used as the mean vector of a multivariate Bernoulli distribution from which 2,000 samples are drawn and evaluated. The latter indicates a potential advantage of the FW iterates.} 
 
The average running times for solving \eqref{eq:MLEr} and \eqref{eq:MLErc} reported in Table~\ref{tbl:runtime} indicate that the ML verification scheme using the simplified model is roughly seven to ten times faster. All algorithms here were run till convergence and without the early stop procedure described earlier. 

\begin{table}[t]
\renewcommand{\arraystretch}{1.1}
\caption{Average Running Times and Step Sizes \textcolor{black}{for $T=50$}}\vspace*{-0.5em}
\label{tbl:runtime} \centering
\begin{tabular}{|l|rr|rr|}
\hline
\hline
 & \multicolumn{2}{c|}{IEEE 37-bus} & \multicolumn{2}{c|}{IEEE 123-bus}\\ \hline
ML and MAP solvers  & Time [sec] & $\mu$ & Time [sec] & $\mu$\\
\hline\hline
PGD solver for \eqref{eq:MLEr} & 60 &0.007  & 336 & 0.0005  \\
\hline
FW solver for \eqref{eq:MLErc} & 5 &  & 45 & \\
\hline
PGD solver for \eqref{eq:MAP} & 204 & 0.08  & 1332 & 0.001\\
\hline
\hline
\end{tabular}
\end{table}

\textcolor{black}{Our algorithms were compared with a random selection of an admissible topology, and an adapted version of the scheme of~\cite{BoSch13}. The latter computes a pseudo-inverse of the sample voltage covariance. The grid topology is identified by passing the signs of the entries in the latter matrix to a spanning tree algorithm. Since this is an identification scheme and to make the comparison more fair, we checked only existing lines (energized or inactive). The scheme was tested only for the 37-bus feeder, because it cannot not handle disjoint feeders.} 

\begin{figure}[t]
\centering	
\includegraphics[width=0.46\textwidth]{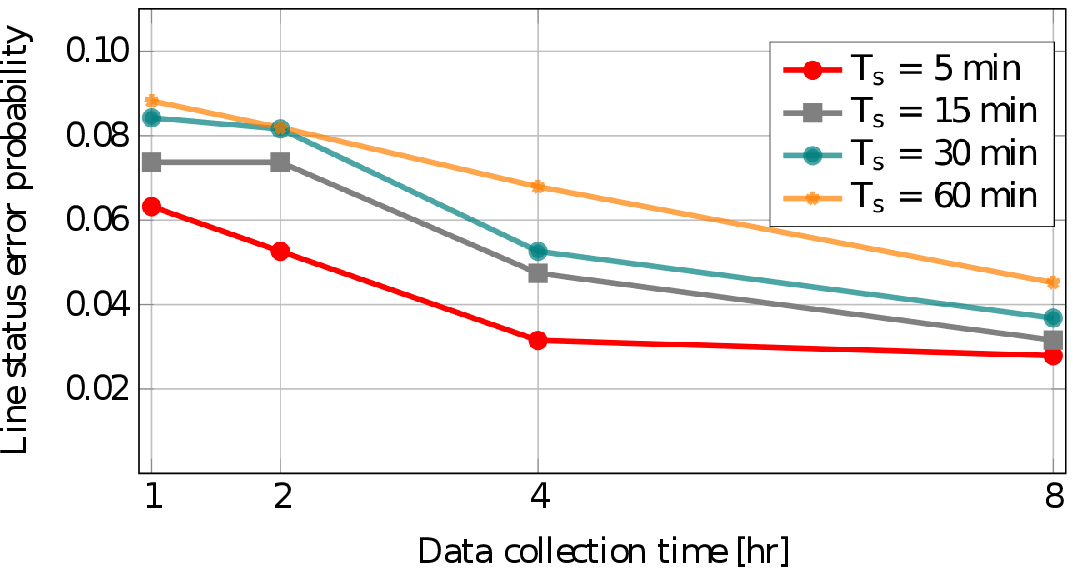}\\
\vspace*{.5em}
\includegraphics[width=0.46\textwidth]{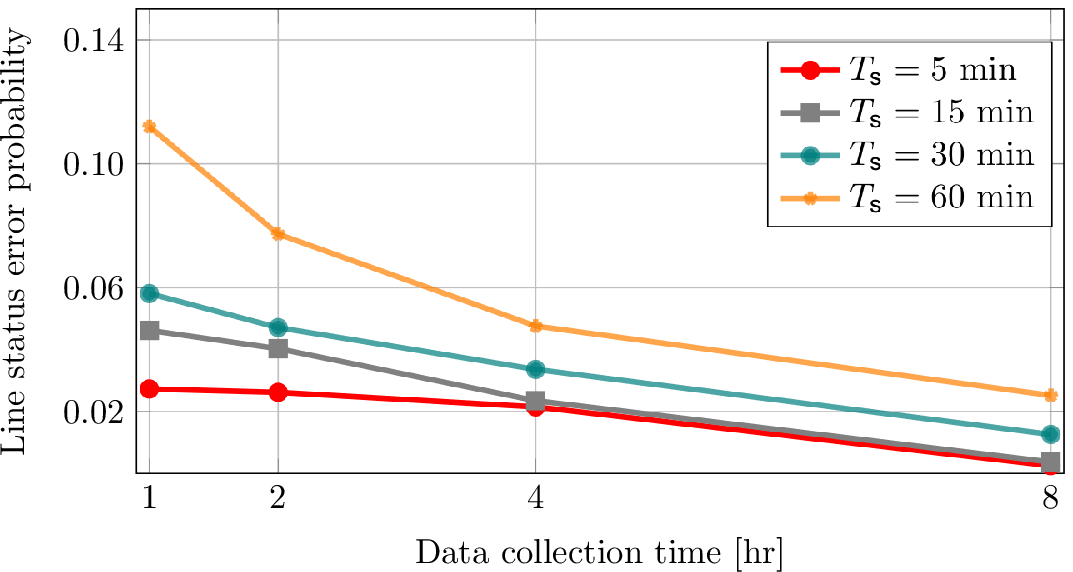}
\caption{Probability of line status errors for the IEEE 37-bus (top) and the IEEE 123-bus (bottom) feeders for varying sampling rates $T_s$.}
\label{fig:res_fixedtimes}
\end{figure}

The average number of line errors for the two feeders is depicted in Fig.~\ref{fig:res}. It is worth emphasizing that reliable line verification can be obtained even for $T<N$, when the sample covariance $\hbSigma$ is \emph{singular}. \textcolor{black}{The performance of the scheme from \cite{BoSch13} is worse than that of a randomly selected topology, since the covariance matrix is away from its ensemble value and line impedances have not been exploited. Figure~\ref{fig:res_fixedtimes} tests the effect of the data sampling rate on the line status error probability achieved by the non-convex ML task of \eqref{eq:MLEr}. The performance improves as $T_s$ and the total collection time increase.}

\begin{figure}[t]
\centering	
\includegraphics[width=0.3\textwidth]{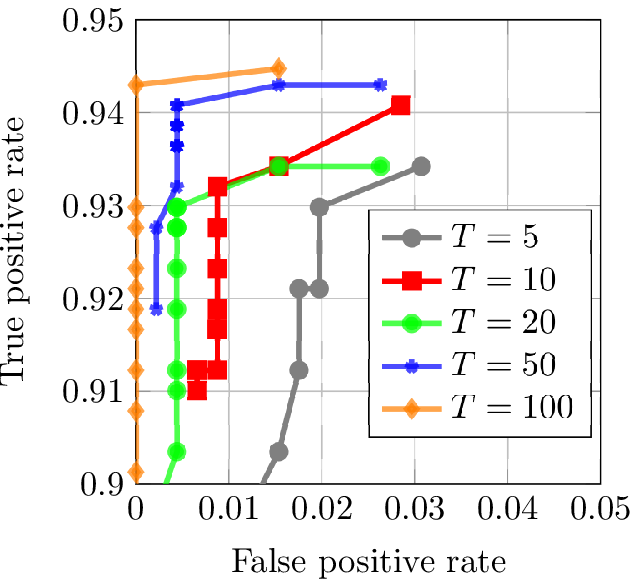}\\
\vspace*{1em}
\includegraphics[width=0.3\textwidth]{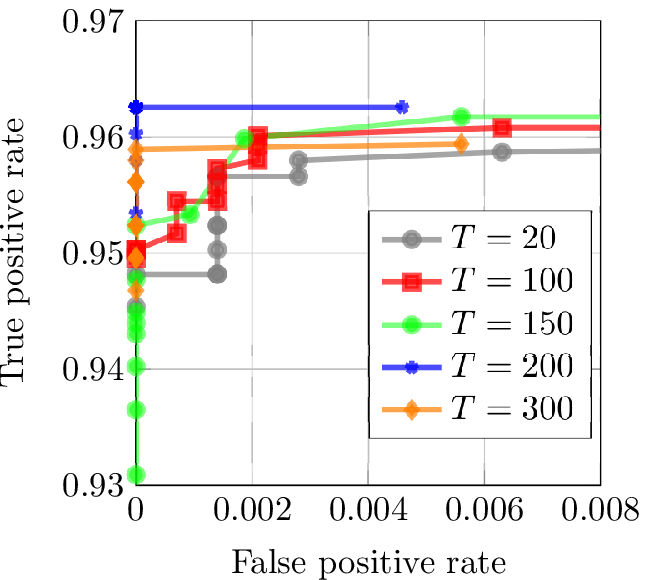}
\caption{\textcolor{black}{True- versus false-positive rates} using varying truncation thresholds for the IEEE 37-bus (top) and 123-bus (bottom) grids with $T_s = 5$ min.}
\label{fig:roc}
\end{figure}

The MAP approach of Sec.~\ref{subsec:priors} was finally tested via 50 Monte Carlo runs. The prior probabilities $\pi_\ell$ were set to 0.5 for switches and 0.9 for lines. A stationary point of \eqref{eq:MAP} was found through the PGD iterates of \eqref{eq:PGD3}--\eqref{eq:PGD4} after random initialization within $[\mathbf{0},\mathbf{1}]$. The entries of the PGD solution were truncated to binary upon thresholding. Varying the threshold from 0.05 to 0.60 in increments of 0.05 resulted in detection curves of Fig.~\ref{fig:roc}. The threshold yielding the point closer to the top-left corner of these plots was used in the MAP curves of Fig.~\ref{fig:res} and \textcolor{black}{the related running times and step sizes are listed in Table~\ref{tbl:runtime} for $T=50$.} A direct comparison between the ML and MAP estimates may not be fair since the latter uses prior information on lines. 

\section{Conclusions}\label{sec:conclusions}
Distribution grid topology verification has been posed as a statistical inference problem. By judiciously expressing the observed voltage magnitudes as functions of the underlying grid topology, maximum likelihood and maximum \emph{a posteriori} probability detection schemes have been proposed. The novel learning tasks minimize (non)-convex objectives depending on the accuracy of the adopted grid data model. Either way, the optimization is over the non-convex feasible set of active line configurations that is relaxed to its convex hull. Solvers of complementary strengths have been devised. Numerical tests using real data on benchmark feeders demonstrate that our schemes perform well even when the number of data is smaller than the network size. Depending on prior information and load activity, collecting smart meter data over 0.5--3 hours can verify topologies over a hundred of nodes. Extending the approach to multiphase grids and dynamically selecting meters constitute interesting research directions. Coping with the task of topology identification where line parameters are unknown too is also technically challenging and practically relevant. \textcolor{black}{To improve scalability, second-order (Newton) and/or accelerated optimization techniques could be pursued. Since data arrive sequentially, the iterates could be initialized to their most recent value.}

\appendix\label{sec:appendix}
\begin{IEEEproof}[Proof of Lemma~\ref{le:gradient}]
The partial derivatives $\frac{\partial f(\bb)}{\partial b_{\ell}}$ for $\ell\in\mcL_e$ will be derived using the matrix differentiation rules~\cite{HornJohnson}:\\
\hspace*{1em}{(r1)} $\partial (\bX \bY)=(\partial \bX)\bY + \bX(\partial \bY)$;\\
\hspace*{1em}{(r2)} $\partial(\log |\bX|)=\trace(\bX^{-1}\partial \bX)$; and\\
\hspace*{1em}{(r3)} $\partial(\bX^{-1})=-\bX^{-1}(\partial \bX) \bX^{-1}$.\\
Rules (r2) and (r3) together with the commutativity property of the trace provide that
\begin{equation}\label{eq:diff1}
\frac{\partial f(\bb)}{\partial b_{\ell}}=\trace\left( \bF(\bb)\frac{\partial \bSigma(\bb)}{\partial b_{\ell}}\right)
\end{equation}
where $\bF(\bb):=\bSigma^{-1}(\bb)-\bSigma^{-1}(\bb) \hbSigma \bSigma^{-1}(\bb)$.

From the definition of $\bSigma(\bb)$ in \eqref{eq:C} and (r1), it follows
{\color{black}
\begin{align}\label{eq:diff2}
\frac{\partial \bSigma(\bb)}{\partial b_{\ell}}& =\frac{\partial \bR(\bb)}{\partial b_{\ell}}\bSigma_p \bR(\bb) +\bR(\bb)\bSigma_p \frac{\partial \bR(\bb)}{\partial b_{\ell}}\nonumber\\
&~~~+ \frac{\partial \bX(\bb)}{\partial b_{\ell}}\bSigma_q\bX(\bb) +\bX(\bb)\bSigma_q \frac{\partial \bX(\bb)}{\partial b_{\ell}}\nonumber\\
&~~~+ \frac{\partial \bR(\bb)}{\partial b_{\ell}}\bSigma_{pq} \bX(\bb) +\bR(\bb)\bSigma_{pq} \frac{\partial \bX(\bb)}{\partial b_{\ell}} \nonumber\\
&~~~+ \frac{\partial \bX(\bb)}{\partial b_{\ell}}\bSigma_{pq}^\top\bR(\bb) +\bX(\bb)\bSigma_{pq}^\top \frac{\partial \bR(\bb)}{\partial b_{\ell}}.
\end{align}}
Using again rule (r3) for the matrix inverse differential yields
\begin{equation}\label{eq:diff3}
\frac{\partial \bR(\bb)}{\partial b_{\ell}} = -\frac{1}{2r_\ell}\bR(\bb)\ba_\ell\ba_\ell^\top \bR(\bb).
\end{equation}
A similar expression holds for $\frac{\partial \bX(\bb)}{\partial b_{\ell}}$. Plugging \eqref{eq:diff2}--\eqref{eq:diff3} into \eqref{eq:diff1} and exploiting the symmetry of the matrices and the commutativity of the trace yield the result.
\end{IEEEproof}

\begin{IEEEproof}[Proof of Lemma~\ref{le:stationary}]
It can be easily verified that $\bF(\bb_o)=\mathbf{0}$ when $\hbSigma=\bSigma(\bb_o)$, and therefore, $\nabla f(\bb_o)=\mathbf{0}$ from \eqref{eq:gradient}. Hence, the vector $\bb_o$ is a stationary point of $f(\bb_o)$.

It is next proved that $\bb_o$ is also a global minimizer for problems \eqref{eq:MLE} and \eqref{eq:MLEr}. \textcolor{black}{It can be shown that the Kullback-Leibler divergence between two Gaussian probability density functions $p_i(\bx)=\mathcal{N}(\mathbf{0},\bSigma_i)$ for $i=1,2$, with $\bx\in\mathbb{R}^N$ is ~\cite[Ex.~2.13]{Bishop}
\begin{align*}
D(p_1||p_2)&:=\int p_1(\bx)\log \left[\frac{p_1(\bx)}{p_2(\bx)} \right] d\bx\\
&=\frac{1}{2}\left[\log \frac{|\bSigma_2|}{|\bSigma_1|}-N+\trace(\bSigma_2^{-1}\bSigma_1)\right].
\end{align*}
Recall also that $D(p_1||p_2)\geq 0$. Since $\bSigma(\bb)\succ \mathbf{0}$ for all $\bb\in \mathcal{B}_r$ attaining a finite value, matrix $\bSigma(\bb)$ qualifies as a covariance matrix. Selecting $\bSigma_2$ to take the parametric form $\bSigma(\bb)$ and fixing $\bSigma_1$ to the given matrix $\bSigma(\bb_o)\succ \mathbf{0}$, it follows that the objective in \eqref{eq:MLE} is an affine transformation of $D(p_1||p_2)$ as $f(\bb)=2D(p_1||p_2)+\log |\bSigma(\bb_o)|+N$. Since $D(p_1||p_2)$ attains zero only if $\bSigma_1=\bSigma_2$, the minimum value of $f(\bb)$ over all $\bb\in \mathcal{B}_r$ is hence achieved at $\bb=\bb_o$.}
\end{IEEEproof}

\begin{IEEEproof}[Proof of Lemma~\ref{le:tfconvex}]
By definition of $\bX(\bb)$, the matrix $\bX^{-1}(\bb)=\sum_{\ell\in\mcL_e}\frac{b_\ell}{x_\ell}\ba_\ell\ba_\ell^\top$ is symmetric positive semi-definite and depends linearly on $\bb\in\mathbb{R}^L_+$. Although there exist $\bb$ for which the matrix $\sum_{\ell\in\mcL_e}\frac{b_\ell}{x_\ell}\ba_\ell\ba_\ell^\top$ becomes singular, these points are practically avoided while minimizing $\tilde{f}(\bb)$ since the term $-2\log |\bX^{-1}(\bb)|$ tends to infinity.

The negative log-determinant of a positive definite matrix is known to be a strictly convex function~\cite{BoVa04}. The convexity of the second term of $\tilde{f}(\bb)$ stems from the equivalence:
\begin{align}\label{eq:trace2}
\trace(\bY \bSigma_\alpha^{-1} \bY\hbSigma)=\min_{\bZ}&~\trace(\bZ)\\
\text{s.t.}~&\left[\begin{array}{cc}
\bSigma_\alpha^{1/2}\bZ\bSigma_\alpha^{1/2} & \bY \hbSigma^{1/2}\\
\hbSigma^{1/2}\bY & \mathbf{I}
\end{array}\right]\succeq \mathbf{0}\nonumber
\end{align}
where $\hbSigma^{1/2}$ is the square root of $\hbSigma$ and likewise for $\bSigma_\alpha^{1/2}$. Recall that although $\bSigma_\alpha$ is invertible, the sample covariance $\hbSigma$ may be singular, for example when $T<N$.

To prove the equivalence in \eqref{eq:trace2}, heed that by Schur's complement the constraint yields $\bSigma_\alpha^{1/2}\bZ\bSigma_\alpha^{1/2}\succeq \bY \hbSigma \bY$, or equivalently, $\bZ\succeq \bSigma_\alpha^{-1/2}\bY \hbSigma \bY\bSigma_\alpha^{-1/2}$. Any feasible $\bZ$ can be expressed as $\bZ=\bSigma_\alpha^{-1/2}\bY \hbSigma \bY\bSigma_\alpha^{-1/2} + \bW$ for some $\bW$. Because $\bZ\succeq \bSigma_\alpha^{-1/2}\bY \hbSigma \bY\bSigma_\alpha^{-1/2}$, it holds that $\bW\succeq \mathbf{0}$. Then, the trace of $\bZ$ is strictly larger than $\trace(\bY \bSigma_\alpha^{-1} \bY\hbSigma)$ unless $\bW=\mathbf{0}$. The equivalence follows since the minimizer of \eqref{eq:trace2} is $\bZ^\star=\bSigma_\alpha^{-1/2}\bY \hbSigma \bY\bSigma_\alpha^{-1/2}$.

The constraint in \eqref{eq:trace2} is a linear matrix inequality with respect to $(\bY,\bZ)$, and therefore, convex. Function $\trace(\bY \bSigma_\alpha^{-1} \bY\hbSigma)$ is convex in $\bY$ as the partial minimization of a jointly convex problem over $(\bY,\bZ)$.
\end{IEEEproof}

\bibliographystyle{IEEEtran}
\bibliography{myabrv,power}

\begin{IEEEbiography}[{\includegraphics[width=1in,height=1.25in,clip,keepaspectratio]{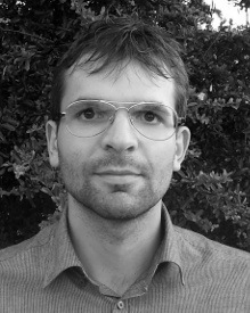}}]{Guido Cavraro} received the B.S. degree in Information Engr., the M.S. degree in Automation Engineering, and the Ph.D. degree in Information Engr. from the Univ. of Padova, Italy, in 2008, 2011, and 2015, respectively. In 2015 and 2016, he was a postdoctoral associate at Information Engr. Dept. of the Univ. of Padova. He is currently a postdoctoral associate at the Bradley Dept. of ECE at Virginia Tech. His current research interests include identification, control and optimization applied to power systems and smart grids.
\end{IEEEbiography} 

\begin{IEEEbiography}[{\includegraphics[width=1in,height=1.25in,clip,keepaspectratio]{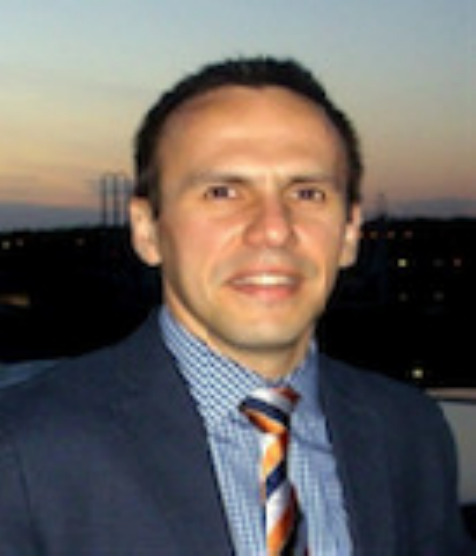}}] {Vassilis Kekatos} (SM'16) is an Assistant Professor of ECE at the Bradley Dept. of ECE at Virginia Tech. He obtained his Diploma, M.Sc., and Ph.D. in computer science and engr. from the Univ. of Patras, Greece, in 2001, 2003, and 2007, respectively. He was a recipient of the Marie Curie Fellowship during 2009-2012, and a research associate with the ECE Dept. at the Univ. of Minnesota, where he received the postdoctoral career development award (honorable mention). During 2014, he stayed with the Univ. of Texas at Austin and the Ohio State Univ. as a visiting researcher. His research focus is on optimization and learning for future energy systems. He is currently serving in the editorial board of the IEEE Trans. on Smart Grid.
\end{IEEEbiography}

\begin{IEEEbiography}[{\includegraphics[width=1in,height=1.25in,clip,keepaspectratio]{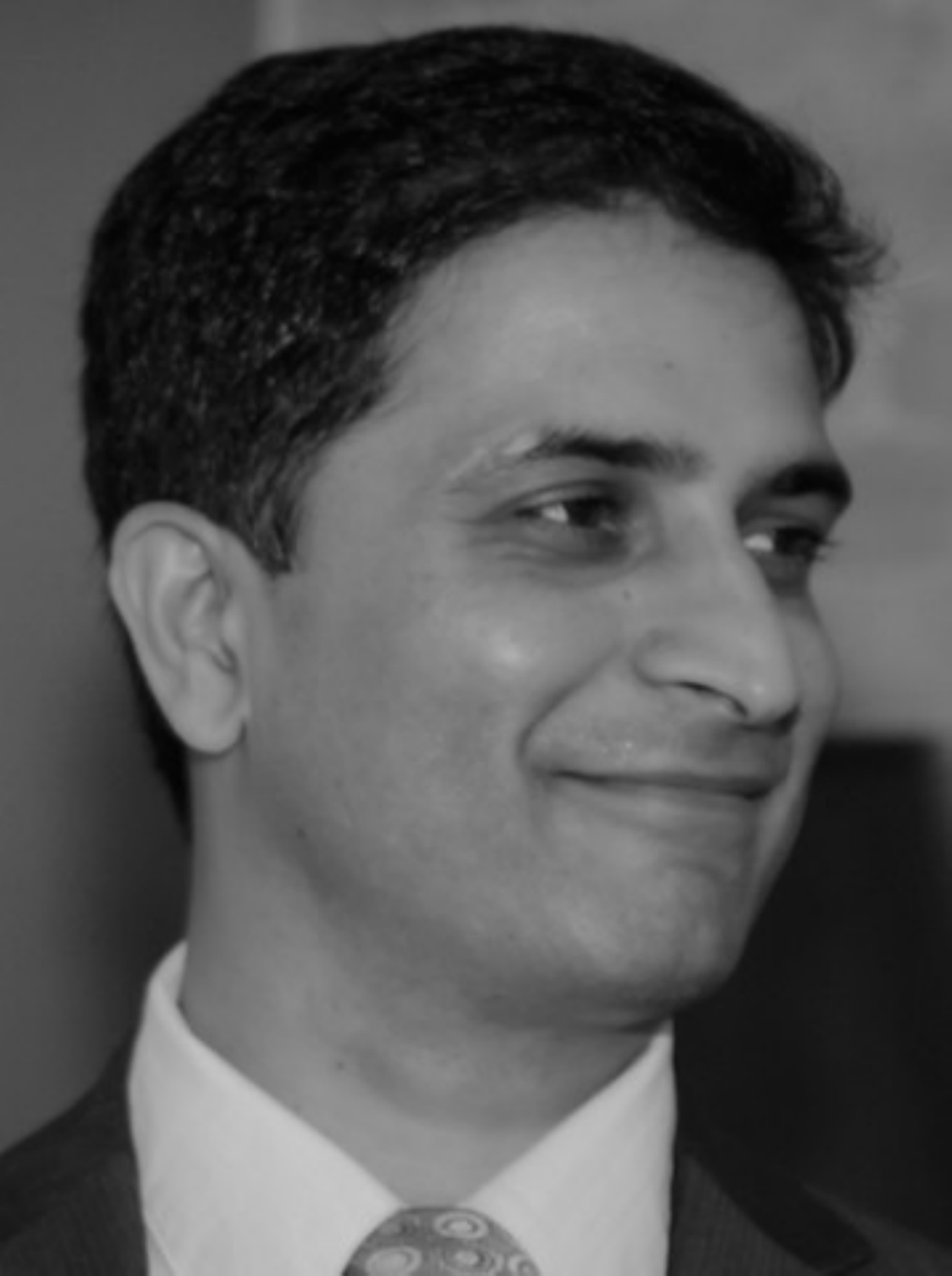}}] {Sriharsha Veeramachaneni} is a researcher at Windlogics Inc., a subsidiary of NextEra Energy Inc., specializing in statistics, machine learning and optimization. He holds a Ph.D. in computer engineering from the Rensselaer Polytechnic Institute.
\end{IEEEbiography}

\end{document}